# Mathematical reflections on modified fractional counting

# Leo Egghe [(1)] and Ronald Rousseau [(2,3)]

[(1)] Hasselt University, Hasselt, Belgium

E-mail: leo.egghe@uhasselt.be

ORCID: 0000-0001-8419-2932

[(2)] University of Antwerp, Faculty of Social Sciences,

Middelheimlaan 1, 2020 Antwerp, Belgium

E-mail: ronald.rousseau@uantwerpen.be

&

[(3)] KU Leuven, MSI, Facultair Onderzoekscentrum ECOOM,

Naamsestraat 61, 3000 Leuven, Belgium

E-mail: ronald.rousseau@kuleuven.be

ORCID: 0000-0002-3252-2538

## Abstract

We make precise what is meant by stating that modified fractional counting (MFC) lies between full counting and complete-normalized fractional counting by proving that for individuals, the MFC-values are weighted geometric averages of these two extremes.

There are two essentially different ways to consider the production of institutes in multi-institutional articles, namely participation and actual number of contributions. Starting from an idea published by Sivertsen, Rousseau and Zhang in 2019 we present three formulae for measuring the production of institutes in multi-institutional articles. It is shown that the one proposed by Sivertsen, Rousseau and Zhang is situated between the two other ways. Less obvious properties of MFC are proven using the majorization order.

## 1. Introduction

Building on Sivertsen's original concept (2016), the authors of Sivertsen et al. (2019) further developed this idea and proposed a method, called modified fractional counting (MFC), which is used at the aggregate level of universities or countries, and aims at obtaining a balanced representation of research efforts across fields and publication patterns (Sivertsen et al., 2025). On the author level, the idea of MFC was validated (at least to some extent) by Sivertsen et al. (2022). Recently, Donner (2024) raised some important questions about MFC as a methodology, a research evaluation tool and its potential for science policy. Donner's objections were answered in (Sivertsen et al., 2025). In this note, we take Donner's considerations into account, but mainly look at MFC from a purely mathematical point of view, largely ignoring practical issues. Concretely, we only study occurrences in the bylines of publications and do not try to include relative research activity, contrary to the underlying idea of (Sivertsen, 2016) and its elaborations. Hence, our work should not be considered as a reply to Donner, but as a contribution to the mathematical study of MFC as a topic in bibliometrics.

On the level of an individual scientist, the main technical-logical contribution of Sivertsen et al. (2019) is that it is shown that the two best-known methods of giving equal credit to all co-authors of a publication, namely full count (each author receives a full credit) and complete-normalized fractional count (each author of an N-author publication receives a credit equal to 1/N) are just the two extremes of a continuum of author credit methods. We make this statement more precise in the first part of this note. In the second part, we offer some mathematical thoughts on MFC as a bibliometric indicator for institutional production and compare it to two related parametric sequences of indicators.

## 2. One author in a multi-authored publication

Consider an article with N authors, where N is a natural number, at least equal to 1. Given k ≥ 1, it is proposed in Sivertsen et al. (2019), to give an author in such an article a credit of

$$MFC_k = \frac{1}{N^{1/k}} = \left(\frac{1}{N}\right)^{1/k} \qquad (1)$$

Clearly, MFC is a parameter family of indicators, with parameter k. Note that k does not have to be an integer. If k = 1 we obtain the case of complete-normalized fractional counting, and in the limit for $k \to \infty$, one obtains the case of full or total counting. We will denote this case simply as $MFC_\infty$. It is a well-known mathematical fact that

$$\forall\, k \geq 1:\ \frac{1}{N} = MFC_1 \leq \left(\frac{1}{N}\right)^{\frac{1}{k}} = MFC_k \leq 1 = MFC_\infty \qquad (2)$$

In order to show how exactly $MFC_k$ fits into a continuum of author credit methods, we recall the following definitions.

Definitions: weighted averages (Hardy et al., 1934, 2.2)

Let $X = (x_1, \ldots, x_n)$ be a vector in $(\mathbb{R}_0^+)^n$, where $\mathbb{R}_0^+$ denotes the strict positive real numbers, and let $W = (w_1, \ldots, w_n)$ be any vector in $[0,1]^n$,

such that $\sum_{j=1}^{n} w_j > 0$. The vector W will be used as a weight vector. Then we define:

the weighted arithmetic average A(X) as:

$$A(X) = \frac{\sum_{j=1}^{n} w_j x_j}{\sum_{j=1}^{n} w_j} \qquad (3)$$

the weighted geometric average G(X) as:

$$G(X) = \left(\prod_{j=1}^{n} x_j^{w_j}\right)^{1/(\sum_{j=1}^{n} w_j)} \qquad (4)$$

and the weighted harmonic average H(X) as:

$$H(X) = \frac{\sum_{j=1}^{n} w_j}{\sum_{j=1}^{n} \left(\frac{w_j}{x_j}\right)} \qquad (5)$$

We will use these definitions for the case $n = 2, X = (MFC_1, MFC_\infty)$ and $W = (\lambda, 1 - \lambda)$, with $0 \leq \lambda \leq 1$. The definitions for averages (3), (4), and (5) then become:

$$A_\lambda(X) = \lambda\, MFC_1 + (1 - \lambda)\, MFC_\infty \qquad (6)$$

$$G_\lambda(X) = (MFC_1)^\lambda \cdot (MFC_\infty)^{1-\lambda} \qquad (7)$$

$$H_\lambda(X) = \frac{1}{\frac{\lambda}{MFC_1} + \frac{1-\lambda}{MFC_\infty}} \qquad (8)$$

Next, we present the following remarkable result.

Theorem 1. $\forall \lambda \in [0,1]: G_\lambda = MFC_{(1/\lambda)}$

Proof. $\forall \lambda \in [0,1]: G_\lambda = (MFC_1)^\lambda \cdot (MFC_\infty)^{1-\lambda} = \left(\frac{1}{N}\right)^\lambda = MFC_{(1/\lambda)}$ □

This theorem shows that every $MFC_k$ is a geometric average of the two extremes with k = 1/λ, and this independent of N.

Corollary

$MFC_2$ (the square root case) corresponds to using equal weights (λ = ½) in the classical (unweighted) geometric average of Theorem 1.

Remarks. We will next show that a similar result involving the arithmetic or harmonic average is not possible.

Proposition 1

It is not possible to write any $MFC_k$ (1 < k < ∞) as a (weighted) arithmetic average of $MFC_1$ and $MFC_\infty$, independent of N.

Proof. Assume that $A_\lambda = \lambda. MFC_1 + (1-\lambda). MFC_\infty = \frac{\lambda}{N} + (1-\lambda)$. Requiring this expression to be equal to some $MFC_k$, we obtain:

$$\left(\frac{1}{N}\right)^{1/k} = \frac{\lambda}{N} + (1-\lambda) \qquad (9)$$

leading to: $k = \frac{ln\left(\frac{1}{N}\right)}{ln\left(\frac{\lambda}{N}+(1-\lambda)\right)}$ . This expression depends on N (and on $\lambda$). □

Proposition 2.

It is not possible to write any $MFC_k$ (1 < k < ∞) as a (weighted) harmonic average of $MFC_1$ and $MFC_\infty$, independent of N.

Proof. Assume that $H_\lambda = \frac{1}{\frac{\lambda}{MFC_1}+\frac{1-\lambda}{MFC_\infty}} = \frac{1}{\lambda N + 1 - \lambda}$. Requiring this expression to be equal to some $MFC_k$, we need:

$$\left(\frac{1}{N}\right)^{1/k} = \frac{1}{\lambda N + 1 - \lambda} \qquad (10)$$

and hence: $k = \frac{\ln(N)}{\ln(\lambda N + 1 - \lambda)}$. Again, we have an expression that depends on N (and on $\lambda$).□

The sum of the individual credits of an N-author article is

$$N \cdot \frac{1}{\sqrt[k]{N}} = N^{((k-1)/k)} \qquad (11)$$

In particular, if k=1 this sum is 1 and when k tends to infinity, the sum tends to N. The following Table 1 provides some values for different N and different k.

Table 1. Author credits and their sums for an N-author article.

| number of authors | N | 1 | 2 | 3 | 5 | 10 | 100 |
|---|---|---|---|---|---|---|---|
| credit of one author | k=1 | 1.00 | 0.50 | 0.33 | 0.20 | 0.10 | 0.01 |
| | 2 | 1.00 | 0.71 | 0.58 | 0.45 | 0.32 | 0.10 |
| | 3 | 1.00 | 0.79 | 0.69 | 0.58 | 0.46 | 0.22 |
| | 5 | 1.00 | 0.87 | 0.80 | 0.72 | 0.63 | 0.40 |
| | 10 | 1.00 | 0.93 | 0.90 | 0.85 | 0.79 | 0.63 |
| | ∞ | 1.00 | 1.00 | 1.00 | 1.00 | 1.00 | 1.00 |
| | | | | | | | |
| Total credit of one publication | k=1 | 1.00 | 1.00 | 1.00 | 1.00 | 1.00 | 1.00 |
| | 2 | 1.00 | 1.41 | 1.73 | 2.24 | 3.16 | 10.00 |
| | 3 | 1.00 | 1.59 | 2.08 | 2.92 | 4.64 | 21.54 |
| | 5 | 1.00 | 1.74 | 2.41 | 3.62 | 6.31 | 39.81 |
| | 10 | 1.00 | 1.87 | 2.69 | 4.26 | 7.94 | 63.10 |
| | ∞ | 1.00 | 2.00 | 3.00 | 5.00 | 10.00 | 100.00 |

We note that, as stated in (Sivertsen et al., 2019, 2025), the sum of all credits of an article is not considered an element in the practical application of MFC for aggregate units. Yet, mathematically, this sum can always be calculated, and we will show that sums do play a role in practice, see section 5 on majorization.

Figures 1 and 2 show $MFC_k = \left(\frac{1}{N}\right)^{1/k}$ and $G_\lambda = MFC_{1/\lambda} = \left(\frac{1}{N}\right)^\lambda$ for the case N=10.

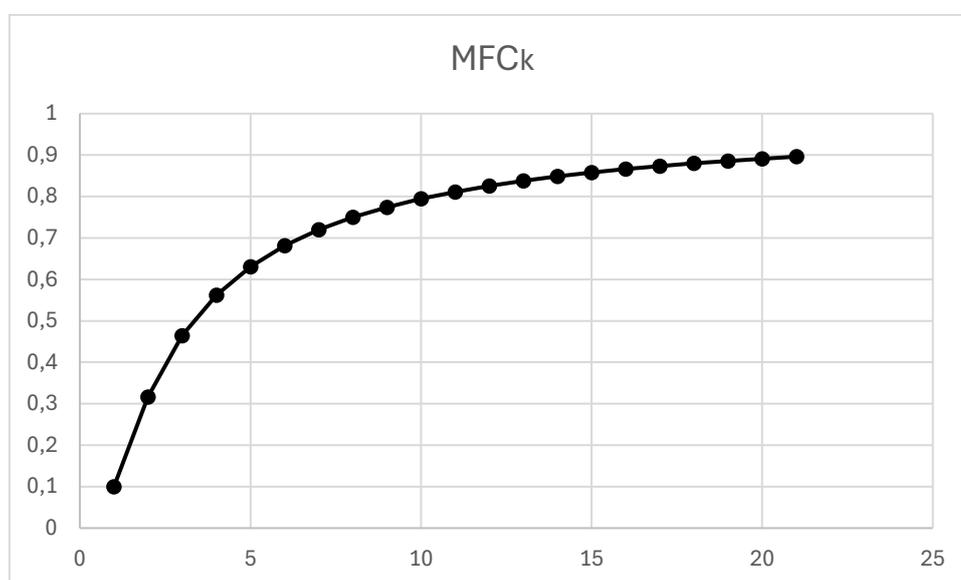

Fig. 1. For any value of N > 1, this curve is concavely increasing

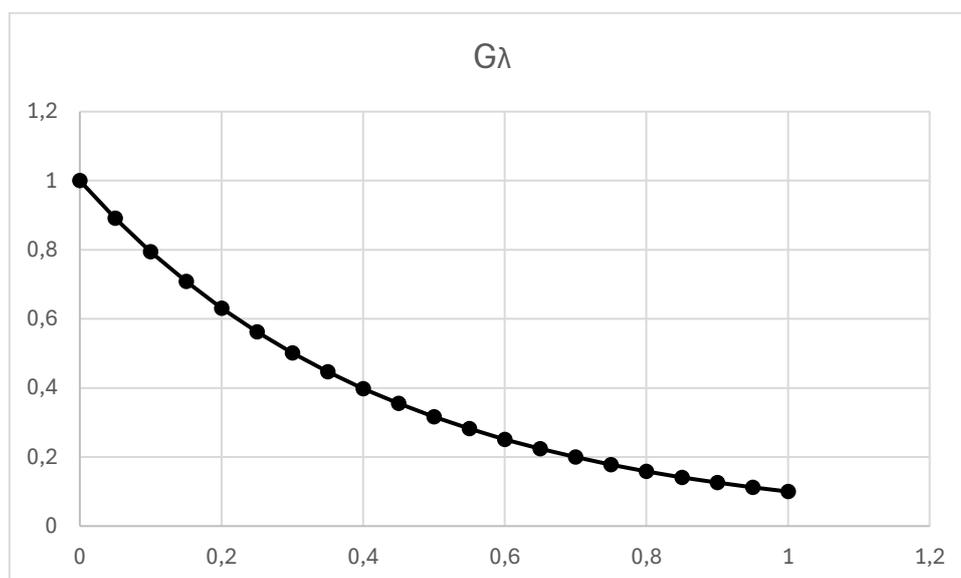

Fig. 2. For any value of N, this curve is convexly decreasing

## 3. Contributions of institutes to a set of multi-institutional articles

Before elaborating on this case formally, we give some informal comments related to Donner (2024) and the original article by Sivertsen et al. (2019).

*3.1 An informal approach*

We note that there are, classically two essentially different ways of studying the publications of an institute (research group, university, country), namely studying participation in (a yes or no situation) or actual contribution to publications. Traditionally these two approaches can be measured either by fractional counting or by full counting (Gauffriau, 2021). In the case of participations one just establishes if the institute is mentioned in the byline. With full counting, each institute receives a credit of 1 and with fractional counting, it receives a credit of 1/(number of participating institutes) for each publication under study. In the case of contributions, one takes the number of times the institute is mentioned in the byline into account (we ignore here, for simplicity, scientists with multiple affiliations). With full counting, the institute receives a credit equal to the number of authors (mentions) in the byline, and with fractional counting, the institute receives a credit equal to the number of mentions divided by the total number of authors. Following Gauffriau (2021) the score functions in the case of participation are called whole counting and whole-fractionalized counting, and in the case of contributions, they are referred to as complete counting and complete-fractionalized counting.

Consider an example of a publication with 7 authors: 4 from institute A, 2 from institute B and one from institute C. Table 2 gives the scores for institutes A, B, and C, using the four score functions mentioned above.

Table 2. Scores for the four classical counting methods

|  | A | B | C | Total score |
|---|---|---|---|---|
| Complete | 4 | 2 | 1 | 7 |
| Fractionalized-complete | 4/7 | 2/7 | 1/7 | 1 |
| Whole | 1 | 1 | 1 | 3 |
| Fractionalized-whole | 1/3 | 1/3/ | 1/3 | 1 |

If we want to study sequences of indicators, depending on a parameter k, then the values in Table 2 correspond to k = 1 and k = ∞ (details follow). Now, we discuss how to assign values for intermediate values of k. First, we can keep the total score, namely $N.\frac{1}{\sqrt[k]{N}}$, and give 4/7th, 2/7th and 1/7th of this score to each institute. For the example above and k = 2, we have $\frac{4}{7} * 7.\frac{1}{\sqrt{7}} = \frac{4}{\sqrt{7}}; \frac{2}{7} * \sqrt{7} = \frac{2}{\sqrt{7}}; \frac{1}{7} * \sqrt{7} = \frac{1}{\sqrt{7}}$.

Yet, as Paul Donner (2024) noted, this is not what has been done in (Sivertsen et al., 2019). In their example, reexamined by Paul Donner, the authors wrote "if a publication has been co-authored by five researchers and two of these researchers are affiliated with the university in focus, the university is credited with a score of 0.63 with MFC$_2$", while the previous method would yield $\frac{2}{5} * \sqrt{5} = 0.89$, a result also mentioned by Paul Donner.

What was actually done to obtain the value 0.63, was also adapting the weights depending on the parameter k, and hence also the total score of a publication.

According to this logic, the weights 4,2,1 become the weights $\sqrt[k]{4}, \sqrt[k]{2}, \sqrt[k]{1}$ and the scores for each institute become: $\sqrt[k]{\frac{4}{7}}, \sqrt[k]{\frac{2}{7}}, \sqrt[k]{\frac{1}{7}}$. For the case k = 2, this becomes $\sqrt{\frac{4}{7}}, \sqrt{\frac{2}{7}}, \sqrt{\frac{1}{7}}$. In this approach with adapted weights, there is no simple mathematical formula for the total score of a publication. If k = 1 we have scores 4/7, 2/7, 1/7 and in the limit for k to infinity, we have 1,1,1. Applying this formula, the total score of a publication goes from 1 (for k=1) to the number of different institutions participating in the target publication (for k = infinity).

Finally, if a group of publications, such as the set of all publications of a country, is the object of study, the score of that country is simply the sum of the scores garnered by that country in each publication under study.

*3.2 A formal approach*

Let P = {$p_1, …, p_n$} be a set of publications, let S be an institute, and let {$q_1,…,q_m$} be the set of authors working in S that participate in at least one of the publications in P. We simply write this as S = {$q_1,…,q_m$}. Of course, not all authors of publications in P belong to S, otherwise, there was no problem as we would not have a multi-institutional framework. Now we define the (m$_x$n) incidence matrix $M(S) = (a_{ij})$, with i =1, …, m, j = 1,…, n, as

$$a_{ij} = 1 \text{ if } q_i \in p_j \text{ and } a_{ij} = 0 \text{ if } q_i \notin p_j \qquad (12)$$

Next, we consider an example.

Example 1. If P = {$p_1,p_2,p_3$}, S = {$q_1,q_2,q_3,q_4$} and $p_1$= {$q_1,q_3,r_1$}, $p_2$ = {$q_2,r_2$} and $p_3$ = {$q_1,q_3,q_4,r_3,r_4$}, where $r_j$, j =1,…,4 denotes an author not belonging to S. Then M(S) is the following 4$_x$3 matrix:

$$M(S) = \begin{pmatrix} 1 & 0 & 1 \\ 0 & 1 & 0 \\ 1 & 0 & 1 \\ 0 & 0 & 1 \end{pmatrix}$$

This matrix is similar to an incidence matrix in graph theory, but differs, for instance, by the fact that incidence matrices of simple undirected graphs (when columns represent edges) have exactly two ones in each column.

We now use the elements of the matrix *M* to define modified fractional counting for an institute S, offering two alternatives as a starting point.

Notation. Consider publication $p_j$ with $N_j$ authors from $M_j$ different institutes, $Y_j$ of which belong to S (j=1, …, n). Then, by (12)

$$\forall j = 1, \ldots, n\colon Y_j = \sum_{i=1}^{m} a_{ij} \leq N_j \tag{13}$$

Definition

$$CMFC_k(S) = \sum_{j=1}^{n} \left( \sum_{i=1}^{m} \left( \frac{a_{ij}}{N_j} \right)^{1/k} \right) = \sum_{j=1}^{n} \left( \frac{\sum_{i=1}^{m} a_{ij}}{N_j^{1/k}} \right) \tag{14}$$

Here C stands for contributions.

Formula (14) = $\sum_{j=1}^{n} \left( \frac{\sum_{i=1}^{m} a_{ij}}{N_j^{1/k}} \right) = \sum_{j=1}^{n} \left( \frac{\sum_{i=1}^{m} a_{ij}}{N_j} N_j (N_j^{-\frac{1}{k}}) \right) = \sum_{j=1}^{n} \left( \frac{Y_j}{N_j} (N_j^{(1-\frac{1}{k})}) \right)$. This shows that CMFC corresponds to the first method highlighted by Donner (2024) and described here in the "informal" part.

Definition. MFC, an alternative for (14)

$$MFC_k(S) = \sum_{j=1}^{n} \left( \frac{\sum_{i=1}^{m} a_{ij}}{N_j} \right)^{1/k} = \sum_{j=1}^{n} \left( \frac{Y_j}{N_j} \right)^{1/k} \tag{15}$$

Formula (15) is the one used in (Sivertsen et al., 2019).

Next, we investigate how formulae (14) and (15) bridge complete normalized fractional counting and full counting.

Taking k = 1 we see that the two formulae yield the same result, corresponding to the "contributions" case, i.e. we have complete-fractionalized counting. Taking the limit for k tending to infinity, we find the following limits:

$$\lim_{k\to\infty}(CMFC_k(S)) = \sum_{j=1}^{n} \lim_{k\to\infty}\left(\frac{Y_j}{N_j^{1/k}}\right) = \sum_{j=1}^{n} Y_j$$ ; this is the total number of collaborations, corresponding to complete counting. Introducing the notation $\delta_j$ which is equal to 1 if the target institute participates in publication j and is zero otherwise, we find:

$$\lim_{k\to\infty}(MFC_k(S)) = \sum_{j=1}^{n} \lim_{k\to\infty}\left(\frac{Y_j}{N_j}\right)^{1/k} = \sum_{j=1}^{n} \delta_j,$$ which is the total number of participations, corresponding to whole counting.

From these limits and the fact that CMFC and MFC are increasing in k, we see that CMFC bridges fractional and full counting in the collaboration case, or stated otherwise CMFC connects complete-fractional counting to complete counting. MFC connects fractional counting in the collaboration case (complete-fractionalized) with full counting in the participation case (whole counting). This suggests a third formula, bridging fractional and full counting in the participation case, i.e., from whole-fractionalized to whole counting. Denoting this formula as PMFC, we define

$$PMFC_k(S) = \sum_{j=1}^{n}\left(\frac{\delta_j}{M_j}\right)^{(1/k)} \qquad (16)$$

Then, for k =1, we obtain the result for fractionalized-whole counting, namely $\sum_i \frac{1}{M_i}$ where the index i refers to those publications in which institute S participates. For k = ∞, we find

$$\lim_{k \to \infty}(PMFC_k(S)) = \sum_{j=1}^{n} \lim_{k \to \infty} \left(\frac{\delta_j}{M_j}\right)^{1/k} = \sum_{j=1}^{n} \delta_j$$

This is the value for whole counting. We note that $\forall j = 1, \ldots, n$: $\delta_j \leq Y_j = \sum_{i=1}^{m} a_{ij}$, with equality only if $Y_j = 0$ or 1.

Proposition 3.

$$MFC_k(S) \leq CMFC_k(S) \Leftrightarrow k \geq 1 \qquad (17)$$

Proof. Based on formulae (14) and (15) we show for all j = 1, …n and all k:

$$\left(\sum_{i=1}^{m} a_{ij}\right)^{\frac{1}{k}} \leq \sum_{i=1}^{m} (a_{ij})^{\frac{1}{k}} \qquad (18)$$

It is well-known (Hardy et al., 1934) that for an array Z with non-negative values, $Z = (z_i)_{i=1,\ldots,m}$,

$$\sum_{i=1}^{m} z_i \leq \left(\sum_{i=1}^{m} (z_i)^{\frac{1}{k}}\right)^k \text{ iff } k \geq 1$$

which proves this proposition. □

Remark. Inequality (17) is always strict, i.e.,

$$MFC_k(S) < CMFC_k(S) \qquad (19)$$

if k > 1 and there exists at least one publication in P such that S has two or more contributors. This follows from the fact that a number t is equal to $t^{1/k}$, only if t = 1.

An illustration. Consider an institution $S_1$ that collaborates with three members in a 5-author publication; the two other authors belong to the same institution. Then MFC$_k$(S$_1$) = $\sqrt[k]{\frac{3}{5}}$ and CMFC$_k$(S$_1$) = $\frac{3}{\sqrt[k]{5}}$, illustrating the case of strict inequalities in (19). Next, we consider S$_2$. This institution collaborates in three 5-author publications with one member in each. Then MFC$_k$(S$_1$) = $3\sqrt[k]{\frac{1}{5}}$ and CMFC$_k$(S$_2$) is equal to $\frac{1}{\sqrt[k]{5}} + \frac{1}{\sqrt[k]{5}} + \frac{1}{\sqrt[k]{5}} = \frac{3}{\sqrt[k]{5}}$. In this case we have equality, illustrating the special case mentioned after (19). More importantly, although CMFC$_k$(S$_1$) = CMFC$_k$(S$_2$), MFC$_k$(S$_1$) < MFC$_k$(S$_2$). We consider this inequality important, at least from a theoretical point of view, as it illustrates that for MFC outside collaborations are more important than inside collaborations (collaborations within one institute).

Remark 1. If S is a singleton, then PMFC, MFC and CMFC coincide.

Remark 2. There is no fixed inequality between PMFC and MFC. If S has two authors in a 5-author publication, with the other three authors are from the same other institution then: $\forall\, k, 1 \leq k < \infty: PMFC_k(S) = \frac{1}{2^{1/k}} > MFC_k(S) = \sqrt[k]{\frac{2}{5}}$; yet, if S has 2 authors in a 5-author publication, where now the three other authors each belong to a different institute, then $PMFC_k(S) = \frac{1}{4^{1/k}} < MFC_k(S) = \sqrt[k]{\frac{2}{5}}$. Is it possible that $PMFC_k(S) > CMFC_k(S)$? Yes, but not for all k at the same time. Consider a publication with 9 authors, two from institute S and 7 from another institute T. Then $PMFC_k(S) = \frac{1}{2^{1/k}}$ and $CMFC_k(S) = \frac{2}{9^{1/k}}$. For k = 1: $PMFC_1(S) = \frac{1}{2} > CMFC_1(S) = \frac{2}{9}$; and for k = 2 we have $PMFC_2(S) = \frac{1}{\sqrt{2}} \approx 0.707 > CMFC_1(S) = \frac{2}{3} \approx 0.667$. Yet, for k = 3, $PMFC_2(S) = \frac{1}{\sqrt[3]{2}} \approx 0.794 < CMFC_1(S) = \frac{2}{\sqrt[3]{9}} \approx 0.961$

Proposition 4. Replication invariance of MFC

If *(y₁, y₂, …, yᵢ, … , yₘ)* denotes an array of participations of M entities in a publication with N co-authors, and if c > 0, then, for every k, $1 \leq k < +\infty$, and every participating entity $S_{j,\ j\ =\ 1,\ ...,M}$, $MFC_k(S_j) = MFC_k\left(S_i^{(c)}\right)$, where $S_i^{(c)}$ denotes the situation where $y_i$ has become $cy_i$, and hence N has become cN.

Proof. This follows immediately from the definition, formula (15).

In words: the MFC value for an institute that participates with 1 author in a 4-author publication, is the same as that for an institute participating with 3 authors in a 12-author publication, and so on. We think that this is a property worth pursuing, and which is not shared with CMFC.

## 4. Some more properties of MFC

For completeness sake, we first recall the following properties of MFC shown in (Sivertsen et al., 2025).

The first one is on the author level.

Proposition 5 (proposition 1 in (Sivertsen et al., 2025))

If an author is added to an article with N authors, then the score of each original author becomes smaller, with the exception of the case of full counting.

Proposition 6 (Proposition 2 in (Sivertsen et al., 2025))

If the score of an entity (institution, country) is weighted according to the relative number of authors that belong to this entity, and a new entity is added to an article with M entities, then the scores of the original entities decrease.

Proposition 7 (Proposition 3 in (Sivertsen et al., 2025))

If the score of an entity (institution, country) is weighted according to the relative number of authors that belong to this entity, and if an entity wants to add one or more authors, then the scores of all other entities will decrease, again except for the case of full counting.

The practical importance of these propositions is explained in (Sivertsen et al., 2025). Next, we add the following proposition.

Proposition 8.

If the score of an entity (institution, country) is weighted according to the relative number of authors that belong to this entity, and if all entities add the same number of co-authors, then the new score of an entity is larger than or equal to the old one, if and only if the old number of participations of this entity is smaller than or equal than the average number of original participations.

Proof. Let *(a₁, a₂, …, aₘ)* be the original array of participations of M entities, let $T = \sum_{j=1}^{M} a_j$, and let a > 0, be the number added to each $a_j$, then we have to show that

$$\sqrt[k]{\frac{a_j+a}{T+Ma}} \geq \sqrt[k]{\frac{a_j}{T}} \tag{20}$$

$$\Leftrightarrow a_j \leq \frac{T}{M} = \mu \tag{21}$$

We see that equation (20) $\Leftrightarrow \frac{a_j}{T} \leq \frac{a_j+a}{T+Ma} \Leftrightarrow a_j(T + Ma) \leq (a_j + a)T) \Leftrightarrow M a_j \leq T \Leftrightarrow a_j \leq \mu$, which is (21).

Corollary 1. Using the same notation as in Proposition 8, we have $\sqrt[k]{\frac{a_j+a}{T+Ma}} \leq \sqrt[k]{\frac{a_j}{T}} \Leftrightarrow a_j \geq \frac{T}{M} = \mu$

Corollary 2. Unless *(a₁, a₂, …, aₘ)* is a constant array, there exist indices *j₁* and *j₂* such that

$$\sqrt[k]{\frac{a_{j_1}+a}{T+Ma}} < \sqrt[k]{\frac{a_{j_1}}{T}} \text{ and } \sqrt[k]{\frac{a_{j_2}+a}{T+Ma}} > \sqrt[k]{\frac{a_{j_2}}{T}} \quad (22)$$

Proof. This follows immediately from (20), (21), Corollary 1 and the fact that, for a non-constant array *(a₁, a₂, …, aₘ)*,

$$min(a_1.,,,,a_M) < \mu < max(a_1.,,,,a_M) \square$$

This corollary shows that both inequalities can happen as a result of the same action.

Corollary 3. If a ≥ 0, Md = median*(a₁, a₂, …, aₘ)*, and *(a₁, a₂, …, aₘ)* is ranked increasingly, then

$$a_i \leq \frac{Md}{2} \Rightarrow (\forall j; 1 \leq j \leq i): \sqrt[k]{\frac{a_j}{T}} \leq \sqrt[k]{\frac{a_j+a}{T+Ma}} \quad (23)$$

Proof. By Markov's inequality (Chow and Teicher, 1978, p. 85,88) we know that $\frac{Md}{2} \leq \mu$ and hence $\forall j; 1 \leq j \leq i: a_j \leq \mu$. Hence corollary 3 follows immediately from (20),(21).

Note that $a_i \leq \frac{Md}{2}$ may not always happen, as in the case of the array (3,4,4).

## 5. A subtle remark about percentages, explained by majorization

An example. Keep the number of authors fixed (here at 10) as well as the number of authors belonging to the target institute (G) (here at 2). Moreover, also the total number of different institutes in each publication is kept fixed (here at 4). In this way the score of the target institute is constant, but its percentage contribution is not. This is shown in Table 3, for k = 2. Other institutes than G, do not have to be the same in different publications.

Table 3. An example in which only the contributions of other institutes change

| G | | | | MFC$_2$-scores | | | | sum | percentage |
|---|---|---|---|---|---|---|---|---|---|
| 2 | 6 | 1 | 1 | 0.447 | 0.775 | 0.316 | 0.316 | 1.854 | 24.1 |
| 2 | 5 | 2 | 1 | 0.447 | 0.707 | 0.447 | 0.316 | 1.917 | 23.3 |
| 2 | 4 | 3 | 1 | 0.447 | 0.632 | 0.548 | 0.316 | 1.943 | 23.0 |
| 2 | 4 | 2 | 2 | 0.447 | 0.632 | 0.447 | 0.447 | 1.973 | 22.7 |
| 2 | 3 | 3 | 2 | 0.447 | 0.548 | 0.548 | 0.447 | 1.99 | 22.5 |

How can Table 3 be explained? For this, we need the notion of majorization.

Definition: The Lorenz curve (Lorenz, 1905)

Let $X = (x_1, x_2, \ldots, x_N)$ be an N-sequence with $x_j \in \mathbb{R}^+, j = 1, \ldots, N$. If X is an N-sequence, ranked in decreasing order (always used in the sense that ranking is not necessarily strict), then the Lorenz curve of X is the curve in the plane obtained by the line segments connecting the origin (0,0) to the points $\left(\frac{k}{N}, \frac{\sum_{j=1}^{k} x_j}{\sum_{j=1}^{N} x_j}\right)$, k= 1,…,N. For k = N, the endpoint (1,1) is reached.

Definition. The majorization property (Hardy et al., 1934; Marshall et al., 2011).

If X and X' are N-sequences, ranked in decreasing order, then X is majorized by X' (equivalently X' majorizes X), denoted as $X \preccurlyeq_L X'$, if

$$\frac{\sum_{j=1}^{k} x_j}{\sum_{j=1}^{N} x_j} \leq \frac{\sum_{j=1}^{k} x'_j}{\sum_{j=1}^{N} x'_j} \text{ for } k = 1, \ldots, N \tag{24}$$

The index L in $X \leqslant_L X'$ refers to the fact that this order relation corresponds to the order relation between the corresponding Lorenz curves. If $\sum_{j=1}^{N} x_j = \sum_{j=1}^{N} x'_j$, then it is not necessary to divide by this sum in (24).

Continuing now with the explanation of Table 3, we first note that the sum of the scores of one publication is obtained as

$$\sum_{j=1}^{M} \sqrt[k]{b_j} \qquad (25)$$

where $0 \leq b_j \leq 1$, denotes the relative contribution of the j-th institute among the M different institutes. As (3, 3, 2, 2) $\leqslant_L$ (4, 2, 2, 2) $\leqslant_L$ (4,3,2,1) $\leqslant_L$ (5,2,2,1) $\leqslant_L$ (6,2,1,1), and the function in (25) is a diversity measure (sum of concave functions), which are known to respect the majorization order (Patil and Taillie, 1982), we see that the more even the data the larger the total score (sum) of a publication.

This shows that when using MFC an institute scores better (in terms of percentage of the total sum of all scores) the more unequal the role of institutes in its publications (all other aspects being equal).

We recall that the issue of majorization has already been discussed briefly in (Sivertsen et al., 2019).

**6. Discussion**

As stated in the introduction, we wrote this article from a mathematical perspective and did not consider the – important – issue of fairness of practical evaluations. Such considerations, especially for the case k = 2, have been provided by Sivertsen (2016) and in Sivertsen et al. (2019, 2025). Practical consequences of counting methods have also been studied recently in (Korytkowski and Kulczycki , 2019).

When considering the defining formula for MFC(S) we see that it contains, the relative contributions of the members of institute S, namely the values $(b_j) = \left(\frac{Y_j}{N_j}\right)$ in the calculation of $MFC_k(S) = \sum_{j=1}^{n}\left(\frac{Y_j}{N_j}\right)^{1/k}$ . We like to make the important remark that if it were possible to determine the exact – intellectual - contribution of an institute to a publication, not just a relative value based on numerical proportions, then the whole mathematical framework of MFC is still valid. As a first step in this direction one could, in those fields where this applies, in the determination of the $Y_j s$ and their sums $N_j$ , count middle authors as 1, the first author as 4, the second as 2, and the corresponding author as 3 (just as an example).

Returning to Example 1, with  S = {$q_1,q_2,q_3,q_4$} and $p_1$= {$q_1,q_3,r_1$}, $p_2$ = {$q_2,r_2$} and $p_3$ = {$q_1,q_3,q_4,r_3,r_4$}, where $r_j$ , j =1,…,4 denotes an author not belonging to S and assuming that the last author is the corresponding author then we have for institute S,  b-values of 6/9, 4/7, and 7/11, see (22), instead of 2/3, ½ and 3/5.

We emphasize that this is just an example, not a proper proposal, to show that the MFC mathematical framework might also be used in other situations in which institutes are weighted.

## 7. Conclusion

In this article we have made precise what was meant by stating that the indicator family $MFC_k$ lies between full counting and complete-normalized fractional counting for individuals. Indeed, it was shown that for one author all $MFC_k$ are weighted geometric averages of these two extremes.

Three sequences of formulae for measuring the contributions of institutes in multi-institutional articles are presented.

It is shown that next to the formula proposed by Sivertsen (2016) and the generalizing sequence as introduced in Sivertsen et al. (2019) other sequences are equally possible. Practical applications will have to show which is the best, or that all three are, in practice, equally admissible. It is shown that more subtle properties of MFC can be proven using the majorization order.

Finally, it is suggested that the MFC framework can not only be used when contributions of institutes are weighted by numbers of participations, but also by actual intellectual contributions, at least if these could be determined.

**Author contributions**

Conceptualization: Leo Egghe and Ronald Rousseau

Formal analysis: Leo Egghe and Ronald Rousseau

Writing – original draft: Ronald Rousseau

Writing – review & editing: Leo Egghe and Ronald Rousseau

**Conflict of interest**

Ronald Rousseau was a member of the editorial board of QSS when this article was written. He is a member of the board of ISSI. Leo Egghe has no competing interests.

**Funding**

No funding has been received.

**Acknowledgment**

The authors thank the reviewers for their useful remarks, leading to a third method of bridging fractionalized and full counting.